\documentclass[12pt,leqno]{amsart}
\usepackage{amsmath,amssymb,amsthm}
\usepackage{latexsym,amscd}
\newcommand{\ra}{{\rightarrow}}

\newcommand{\M}{{\mathrm {M}}}

\newcommand{\cO}{{\mathcal O}}

\newtheorem{theorem}{Theorem}[section]
\newtheorem{definition}[theorem]{Definition}

\newtheorem{proposition}[theorem]{Proposition}

\newtheorem{corollary}[theorem]{Corollary}

\newtheorem{lemma}[theorem]{Lemma}
\newtheorem{remark}[theorem]{Remark}

\newcommand{\Gal}{{\mathrm {Gal }}}

\newcommand{\sk}{\vspace{0.1in}}

\newcommand{\Res}{\operatorname{Res}}

\newcommand{\Pic}{\mbox{Pic}}

\newcommand{\Norm}{\operatorname{Norm}}
\newcommand{\Tr}{\operatorname{Tr}}
\newcommand{\Aut}{\operatorname{Aut }}

\newcommand{\GL}{{\operatorname{GL}}}
\newcommand{\End}{{\rm{End}}}

\newcommand{\Z}{{\mathbb Z}}
\newcommand{\Q}{{\mathbb Q}}
\newcommand{\C}{{\mathbb C}}
\newcommand{\R}{{\mathbb R}}
\newcommand{\F}{{\mathbb F}}

\newcommand{\qbar}{{\bar\Q}}

\newcommand{\NS}{{\operatorname{NS}}}

\newcommand{\cA}{{\mathcal A}}
\newcommand{\cL}{{\mathcal L}}

\newfont{\gotip}{eufb10 at 12pt}

\newcommand{\ga}{\mbox{\gotip a}}

\include{thebibliography}

\begin{document}

\title{Abelian surfaces of  $\GL_2$-type  as Jacobians of curves}
\author[Josep Gonz\'alez]{Josep Gonz\'alez}
\address{Escola Universit\`aria Polit\`ecnica de Vilanova i la Geltr\'u, Av. V\'{\i}ctor
Balaguer s/n, E-08800 Vilanova i la Geltr\'u, Spain}
\email{josepg@mat.upc.es,\,guardia@mat.upc.es,\,vrotger@mat.upc.es }

\thanks{2000 {\em Mathematics Subject Classification.} Primary 11G10, 11G18, 14H40; Secondary 11F11, 14Q05.}

\thanks{The authors were partially supported by DGICYT Grant
BFM2003-06768-C02-02.}

\author{Jordi Gu\`{a}rdia }

\author{Victor Rotger}

\keywords{Principal polarizations, abelian varieties of $\GL_2$-type, Weil restriction, $\Q$-curves, modular abelian surfaces}

\begin{abstract}
We study the set of isomorphism classes of
principal polarizations on abelian varieties of $\GL _2$-type. As
applications of our results, we
construct examples of curves $C$, $C'/\Q $ of genus
two which are nonisomorphic over $\bar \Q $ and share isomorphic unpolarized modular
Jacobian varieties over $\Q $; we also show a method to
obtain genus two curves over $\Q$ whose Jacobian varieties are isomorphic to Weil's restriction of quadratic $\Q$-curves, and present examples.

\end{abstract}

\maketitle

\section{Introduction}
\sk\noindent Let $C/k$ be a smooth and irreducible projective
algebraic curve over a field $k$ and let $(J(C), \Theta )$ be its
principally polarized Jacobian variety. By Torelli's theorem, $C$
is completely determined by $(J(C), \Theta )$ up to isomorphism
over a fixed algebraic closure $\overline{k}$ of $k$. The question
then arises of whether the unpolarized abelian variety $J(C)$
already determines $C$. It turns out that the answer has its roots
in the arithmetic of the ring of endomorphisms $\End (J(C))$ of
$J(C)$. Indeed, in the generic case where $\End (J(C))=\Z $, the
curve $C$ can be recovered from $J(C)$, but as soon as $\End
(J(C))\varsupsetneq \Z $, it may be very well the case that there
exist finitely many pairwise nonisomorphic curves $C_1$, ...,
$C_{\tau }$ such that
$$
J(C_1)\simeq ... \simeq J(C_{\tau })
$$
as unpolarized abelian varieties. This phenomenon was first
observed by Humbert in \cite{Hu}, where he proved the existence of
two nonisomorphic Riemann surfaces $C_1$ and $C_2$ of genus $2$
sharing the same simple Jacobian variety. Humbert's method was
generalized by Lange in \cite{La}. Finally, it has been recently
shown in \cite{Ro} that there exist arbitrarily large sets $\{
C_1, ..., C_{\tau } \}$, $\tau
>>0$, of pairwise nonisomorphic curves of genus $2$ such that
$J(C_i)$ are isomorphic simple abelian surfaces.

The aim of this paper is to yield an arithmetical and effective
approach to this question by considering abelian surfaces $A_f$
defined over $\Q$ attached by Shimura to a newform $f$ in
$S_2(\Gamma_1(N))$. These surfaces are viewed as optimal quotients
of the Jacobian $J_1(N)$ of the modular curve $X_1(N)$ and are
canonically polarized with the polarization $\cL $ over $\Q $
pushed out from the principal polarization $\Theta _{X_1(N)}$ on
$J_1(N)$. We provide an effective criterion for determining all
their principal polarizations defined over $\Q $ and for
determining which of them correspond to the canonical polarization
of the jacobian of a curve.

In section \ref{gl2}, we consider abelian varieties $A$ of
$\GL_2$-type over $k$, i. e. such that their algebra of
endomorphisms $F=\Q \otimes \End _{k}A$ is a number field of
degree $[F:\Q ]=\dim A$. In Theorem \ref{pi}, we determine the set
of isomorphism classes of principal polarizations on $A$ over $k$
in an explicit and computable way, while in Corollary
\ref{criterion} we provide necessary and sufficient conditions for
$A$ to be principally polarizable over the base field of
definition $k$.

Section \ref{sgl2} is devoted to the study of abelian surfaces of $\GL_2$-type
over a number field $k$ and their principal polarizations from an arithmetical point
of view. As we discuss in detail, these questions are crucially
related to the problem of finding nonisomorphic curves
whose Jacobian varieties are pairwise isomorphic as unpolarized abelian
varieties. One of our main results can be rephrased as follows.

\begin{theorem}\label{main}

Let $A $ be a principally polarizable abelian surface over $\Q$ such that $\End _{\Q
}A=R$ is an order in a quadratic field.  Let $R_0$ be the subring of $R$
fixed by complex conjugation. Then

\begin{enumerate}

\item
If either $R_0=\Z $ or
$R^*$ contains a
unit of negative norm, there is a single isomorphism class of
principal polarizations on $A$ over $\Q $.

\item
Otherwise,
there are exactly two isomorphism classes of principal polarizations on
$A$ over $\Q $. More precisely,
either
\begin{itemize}

\item
there exists a curve $C/\Q $ of genus two and an elliptic
curve $C'$ over a quadratic extension $K$ of $\Q $ such that $A$
is both isomorphic to $J(C)$ and to Weil's restriction $ \Res
_{K/\Q }(C')$ of $C'$ as unpolarized abelian surfaces over $\Q $,
\end{itemize}
or
\begin{itemize}
\item
there exist two  curves $C$, $C'/\Q $ of genus two nonisomorphic
over $\Q$  such that $A$ is isomorphic over $\Q$ to their Jacobian
varieties.
\end{itemize}
\end{enumerate}
\end{theorem}

As a consequence, we prove the existence of infinitely many abelian
surfaces $A/\Q $ of $\GL_2$-type which are simultaneously
the Jacobian variety of a curve $C/\Q$ of genus two and
Weil's restriction of an elliptic curve over a quadratic field.

In the last section explicit examples are presented. First, we
summarize well known results about modular abelian varieties
arising from modular newforms; in the rest of the section we show
how our results can be effectively applied to abelian surfaces
arising from newforms of trivial Nebentypus by using the procedure
described in \cite{GoGoGu}, which allows us to obtain an equation
of a genus two curve from a period matrix in a symplectic basis.
In this way, we construct an explicit example of two nonisomorphic
(over $\bar \Q $) curves $C$, $C'/\Q $ of genus two such that
their Jacobian varieties are isomorphic over $\Q $ to an
absolutely simple quotient of the Jacobian of the modular curve
$X_0(65)$. We also provide explicit examples of curves $C$ of
genus two over $\Q$ whose Jacobian variety $J(C)$ is of
$\GL_2$-type and isomorphic to Weil's restriction over $\Q$ of an
elliptic curve $E$ over a quadratic field. Finally, we exhibit
several pairwise nonisomorphic curves $C_i/\Q $ of genus two such
that their respective Jacobian varieties $J(C_i)$ are mutually
isogenous over $\Q $.

Alternative methods to ours were proposed by Howe (\cite{Ho1},
\cite{Ho2}) and Villegas (\cite{Vi}) and explicit equations of
similar phenomena have been given by Howe in \cite{Ho3} and
\cite{Ho4}.

{\em Acknowledgement.} The third author thanks Gabriel Cardona for
some useful comments on an earlier version of the paper.

\section{Abelian varieties of $\GL_2$-type}\label{gl2}

\subsection{The N\'{e}ron-Severi group of an abelian variety}

We begin  by recalling some basic facts on abelian varieties
(cf.\,\cite{Mi}, \cite{Mu}). Let $k$ be a field, $k^{s}$ a
separable closure of $k$ and let $G_k=\Gal (k^{s}/k)$. For any
field extension $K/k$, let $A_K = A \times _k K$ be the abelian
variety obtained by base change of $A$ to Spec $K$. Let $\Pic (A)$
be the group of invertible sheaves on $A$, $\Pic ^0(A_{k^{s}})$ be
the subgroup of $\Pic (A_{k^{s}})$ of invertible sheaves
algebraically equivalent to $0$ and $\Pic ^0(A)=\Pic (A)\cap \Pic
^0(A_{k^{s}})$. The N\'{e}ron-Severi group of $A$ is $\NS (A) = \Pic
(A)/\Pic ^0(A)$. We shall denote $\NS (A_{k^{s}})^{G_k} = H^0(G_k,
\NS (A_{k^{s}}))$ the group of $k$-rational algebraic equivalence
classes of invertible sheaves. Note that in general not all
elements in $\NS (A_{k^{s}})^{G_k}$ are represented by an element
in $\Pic (A)$.

For any $P\in A(k^{s})$, let $\tau _P$ denote the translation by
$P$ map on $A$. Every invertible sheaf $\cL \in \NS
(A_{k^{s}})^{G_k}$ defines a morphism $\varphi _{\cL }: A \ra \hat
A$ over $k$ given by $\varphi _{\cL }(P)=\tau _P^*(\cL )\otimes
\cL ^{-1}$, which is an isogeny if $\cL $ is nondegenerate. A
polarization on $A$ defined over $k$ is the class of algebraic
equivalence of an ample invertible sheaf $\cL \in \NS
(A_{k^{s}})^{G_k}$. Equivalently, a polarization on $A$ over $k$
is an isogeny $\lambda :A \ra \hat A$ defined over $k$ such that
$\lambda \otimes k^{s} = \varphi _{\cL }$ for some ample line
bundle $\cL $ on $A_{k^{s}}$.

A nondegenerate invertible sheaf $\cL $ on $A$ induces an
anti-involution on the algebra of endomorphisms
$$
*:\Q \otimes \End _{k}A\stackrel{\sim }{\ra }\Q \otimes \End
_{k}A, \quad t\mapsto {\varphi _{\cL }}^{-1}\cdot \hat t\cdot
\varphi _{\cL }.
$$
Let $\End _k^sA=\{ \beta \in \End _k A, \beta ^*=\beta \} $ denote
the subgroup of symmetric endomorphisms with respect to $\cL $ and
$\End _{k +}^sA$ be the set of positive symmetric endomorphisms of
$A$.

{}From now on, we assume that $k$ is a subfield of a fixed
algebraic closure $\bar \Q =\bar k$ of $\Q $.

\begin{proposition}
\label{pol}

Let $A/k$ be an abelian variety and let $\cL \in \NS (A_{\bar
k})^{G_k}$ be nondegenerate.

Then, for any endomorphism $t\in \End_k^s A$, there exists a
unique $\cL ^{(t)}\in \NS (A_{\bar k})^{G_k}$ such that $\varphi
_{\cL ^{(t)}}= \varphi _{\cL }\cdot t$. The following properties
are satisfied:

\begin{itemize}

\item [(i)]
Let $E$ and $E_t$ denote the alternating Riemann forms attached to
$\cL$ and $\cL^{(t)}$ respectively. Then
$$
E_t(x,y)=E(x, t y)=E(t x, y)\,.
$$

\item[(ii)]
If $t$ is a totally positive element, then
$\cL$ is a polarization if and only if
$\cL^{(t)}$ is.

\item[(iii)]
For any $t\in \End _k^s A$ and $d\in \Z $,
$t^*(\cL )=\cL ^{(t^2)}$ and $\cL^{(d)}=\cL ^{\otimes d}$.

\end{itemize}

\end{proposition}

{\em Proof. } Let $A(\C )=V/\Lambda $ for a complex vector space
$V$ and a lattice $\Lambda $. Since the involution induced by $\cL
$ on $\End _k^s A$ is the identity map, it may be checked that,
for any endomorphism $t\in \End _k^s A$, $E_t:V \times V
\rightarrow \R $, $(x, y)\mapsto E(t x, y)=E(x, t y)$ is again an
alternating Riemann form on $A$. Then, by the Appell-Humbert
theorem, there exists a unique invertible sheaf $\cL ^{(t)}$ up to
algebraic equivalence such that $E_{\cL ^{(t)}}=E_t$. It follows
from the analytical representation of the morphism $\varphi _{\cL
}$ in terms of $E$ that $\varphi _{\cL ^{(t)}}= \varphi _{\cL
}\circ t$. Since both $t$ and $\cL $ are defined over $k$, the
same holds for $\cL ^{(t)}$.

As for (ii), let us denote by $H$ the hermitian form on $V$
attached to $\cL$. Then, the matrix of the hermitian form $H_t$
attached to $\cL^{(t)}$, with respect to any basis of $V$, is the
product of the matrices of $H$ and $t$. Since $t$ is a
self-adjoint endomorphism on the hermitian space $(V,H)$, there is
an orthogonal basis of $V$ over which $t$ has diagonal form. It is
clear that if $t$ is totally positive, then $H_t$ is positive
definite if and only if $H$ is.

Finally, we know that for any endomorphism $t\in \End _k^s A$ and
any invertible sheaf $\cL $ on $A$, $E_{t^*(\cL )}(x, y) = E (t x,
t y) = E (t^2 x, y) = E_{\cL ^{(t^2)}}(x, y)$ and thus $t^*(\cL
)=\cL ^{(t^2)}$. Since $E_{\cL ^{\otimes d}} = d E$, we also have
that $\cL^{(d)}=\cL ^{\otimes d}$. $\Box $

\begin{remark} Note that our construction of $\cL ^{(t)}$ extends
for arbitrary $t\in \Q\otimes \End_k A$ to produce elements
 $\cL ^{(t)}\in \Q \otimes \NS (A_{\bar
k})^{G_k}$.
\end{remark}

\begin{theorem}
\label{eps}

Assume there  exists  a principal polarization  $\cL _0$  on
$A$ defined over $k$ and let $\End _k^sA$ denote the subgroup of
symmetric endomorphisms with respect to $\cL_0$. Then there is an
isomorphism of groups

$$
\begin{matrix}
\epsilon : &\NS (A_{\bar k})^{G_k}&\stackrel{\sim }{\ra }&\End _k^sA\\
         &   \cL   &\mapsto       &  \varphi _{\cL _0}^{-1}\cdot
          \varphi _{\cL }
\end{matrix}
$$
such that $\cL \in \NS (A_{\bar k})^{G_k}$ is a polarization if
and only if $\epsilon (\cL )\in \End _{k +}^sA$ and it is
principal if and only if $\epsilon (\cL )\in \Aut _{k +}^sA$.
Moreover, $\epsilon^{-1}(t)=\cL_0^{(t)}$.

\end{theorem}

{\em Proof. } The first part of the theorem is well known if we
replace $k$ by $\bar k$. It is clear that if $\cL _0\in \NS
(A_{\bar k})^{G_k}$, $\epsilon (\cL )\in \End _k^s A$ for any $\cL
\in \NS (A_{\bar k})^{G_k}$. From the relation $$\epsilon
(\cL_0^{(t)} ) =\varphi _{\cL _0}^{-1}\cdot \varphi _{\cL_0^{(t)}
}= \varphi _{\cL _0}^{-1}\cdot \varphi _{\cL_0}\cdot t=t\,,$$ we
obtain that the morphism $\epsilon $ is surjective and the
statement is immediate. $\Box $

\begin{definition}
We shall say that two polarizations $\cL _1$, $\cL _2$ defined
over a field extension $K/k$ are $K$-isomorphic, $\cL
_1\stackrel{K}\simeq \cL _2$, if and only if there exists $u\in
\Aut _K A$ such that $u^*(\cL _2)=\cL _1$. We denote by $\Pi(A_k)$
the set of $k$-isomorphism classes of principal polarizations on
$A_k$ defined over $k$, i. e.,
$$
\Pi(A_k) = \{ \mbox{ principal polarizations } \cL \mbox{ defined
over }k \} /\stackrel {k}\simeq .
$$
\end{definition}

\begin{proposition} \cite{Mi} The set $\Pi(A_k)$ is finite.
\end{proposition}

We will denote by $\pi(A_k)$ the cardinality of $\Pi(A_k)$. Two
positive symmetric endomorphisms $\beta _1$, $\beta
_2\in \End _{k +}^sA$ are equivalent, $\beta _1\sim \beta _2$, if
$\beta _1= \beta ^*\beta _2 \beta $ for some $\beta \in \Aut
_{k}A$. The next result follows from Theorem \ref{eps}.

\begin{corollary}
\label{bij}

The morphism $\epsilon $ induces a bijection of finite sets
$$
\Pi(A_k)\longleftrightarrow \Aut _{k +}^sA/\sim \,.
$$

\end{corollary}

\subsection{Polarizations on abelian varieties of $\GL_2$-type}\label{pgl2}

We now describe explicitly the isomorphism of Theorem \ref{eps}
for a certain class of abelian varieties. This will provide a
procedure for determining all isomorphism classes of principal
polarizations on them in a computable way. We will illustrate it
with several examples in section \ref{examples}.

\begin{definition}

An abelian variety $A$  is of $\GL _2$-type over $k$ if $A$ is
defined over $k$ and $\End _k A$ is an order in a  number field
$F$ of degree $[F:\Q ]=\dim A$. If $F$ is totally real, we then
say that $A$ is of real $\GL _2$-type over $k$.

\end{definition}

\begin{remark}

By Albert's classification of involuting division algebras (see
\cite{Mu}), the algebra of endomorphisms $F=\Q \otimes \End _k A$
of an abelian variety $A$ of $\GL _2$-type over $k$ is isomorphic
either to a totally real field or a CM-field. The Rosati
involution with respect to any polarization $\cL $ on $A$ over $k$
acts as complex conjugation on the number field $F=\Q \otimes \End
_k A$ (see \cite{Mu}) and  it can be checked that the same holds
for the involution with respect to a non necessarily ample
invertible sheaf $\cL $. In particular, $\Q \otimes \End _k^s A$
is isomorphic either to $F$ or the maximal totally real subfield
$F_0$ of $F$, depending on whether $A$ is of real $\GL _2$-type or
not, and independently of the choice of $\cL $. This allows us in
this section to refer to $\End^s_k(A)$ without further mention to
the chosen polarization.

\end{remark}

\begin{remark}
An abelian variety $A/k$ is said to have real multiplication by a
field $F$ if $F$ is a totally real field contained in $\Q\otimes
\End _k A$ with $[F:\Q ]=\dim A$. It is clear that the condition
of being of real $\GL_2$-type is stronger than that of being of
real multiplication.
\end{remark}

The main sources of abelian varieties  of $\GL _2$-type over $\Q$
are the abelian varieties $A_f$ attached by Shimura to an
eigenform $f$ of weight $2$ for the modular congruence subgroup
$\Gamma _1(N)$ and the  real case occurs when $f\in
S_2(\Gamma_0(N))$. Actually, Ribet has proved  in \cite{Ri} that,
assuming Serre's conjecture $(3.2.4_?)$ of \cite{Se}, these  are
all the abelian varieties of $\GL _2$-type over $\Q $.

For the rest of this section, we fix the following notation. For any
totally real number field $F_0$ and any subset $S$ of $F_0$, we shall
denote by $S_+$ the set of totally positive elements of $S$.

The next theorem provides an explicit description of all
polarizations on an abelian variety of $\GL_2$-type.

\begin{theorem}
\label{pi}

Let $A$ be an abelian variety of $\GL_2$-type over $k$. Let
$R=\End_k A$, $R_0$ be the subring of $R$ fixed by complex
conjugation, $F=\Q\otimes R$ and $F_0=\Q\otimes R_0$. Suppose that
$A$ admits a principal polarization $\cL _0 $ on $A$ defined over
$k$ and let $\epsilon$ be as in Theorem \ref{eps}. Then,

\begin{itemize}
\item[(i)]
For any $\cL \in \NS (A_{\bar k})^{G_k}$ and any endomorphism
$s\in R$:
$$
\deg (\cL )=
\vert \Norm _{F/\Q } (\epsilon (\cL ))\vert ,
$$
$$
\epsilon (s^*(\cL )) =
  \begin{cases}
    s^2\cdot \epsilon (\cL ) & \text{if } F=F_0, \\
    \Norm _{F/F_0}(s)\cdot \epsilon (\cL ) & \text{if } [F:F_0]=2.
  \end{cases}
$$

\item[(ii)]
The class of an invertible sheaf $\cL \in \NS (A_{\bar k})^{G_k}$
is a polarization if and only if $\epsilon (\cL )\in (R_0)_{+}$.
In particular, $\cL $ is a principal polarization if and only if
$\epsilon (\cL )\in (R_0^*)_+$.

\item[(iii)]
Let $P(R)$ be the multiplicative group
$$
P(R)=\begin{cases}
   (R_0)_+^*/ R_0^{* 2} & \text{if } F=F_0, \\
    (R_0)_+^*/\Norm_{F/F_0}(R^*) & \text{if } [F:F_0]=2.
  \end{cases}
$$
Then, there is a bijective correspondence
between the sets $\Pi(A_k)$ and $P(R)$
and hence
$$
\pi(A_k) = 2^{N} \mbox{ with } 0\leq N\leq [F_0:\Q ]-1.
$$
The $2^N$ isomorphism classes of principal polarizations on $A$
over $k$ are represented by the invertible sheaves $\cL
_0^{(u_j)}$, where $\{ u_j\} _{j=1}^{2^N}$ is a system of
representatives of $P(R)$.
\end{itemize}

\end{theorem}

{\em Proof. } By Theorem \ref{eps}, for every $\cL \in \NS
(A_{\bar k})^{G_k}$, there exists $t\in R_0$ such that $\cL =\cL
_0^{(t)}$. Thus, $\epsilon (s ^*(\cL )) = \epsilon (s^*(\cL
_0^{(t)}))=\epsilon (\cL _0^{( s \bar s t)})= s \bar s \cdot
\epsilon (\cL )$, where $\bar {}$ denotes complex conjugation.
Moreover, $\deg (\cL )=\deg (\epsilon (\cL ))^{1/2}=\vert \Norm
_{F/\Q } (\epsilon (\cL ))\vert $. This yields (i).

Part (ii) follows immediately from Theorem \ref{eps} and part (i).

Finally, the bijection $\Pi(A_k)\simeq  P(R)$ is now a consequence
of Corollary \ref{bij}. Dirichlet's unit theorem then implies that
$\pi (A_k) = 2^{N}$ for $0\leq N\leq [F_0:\Q ]-1$. From (i) we
obtain that any two polarizations $\cL $, $\cL '$ on $A$ defined
over $k$ are isomorphic if and only if $\epsilon (\cL )= u \bar
u\,\epsilon (\cL ')$ for some $u\in R^*$. Hence, representatives
of the isomorphism classes of principal polarizations on $A$ over
$k$ are provided by $\cL _0^{(u_j)}$ for representatives $u_j$ of
$P(R)$. $\Box $

\begin{remark}\label{degree1}

If $\cL_0 \in \NS (A_{\bar k})^{G_k}$ is a nonnecessarily ample
invertible sheaf on $A$ over $k$ of degree $1$, then the map
$\epsilon : \NS (A_{\bar k})^{G_k}\rightarrow \End _kA$, defined
by $\epsilon (\cL)= \varphi_{\cL_0}^{-1}\cdot \varphi (\cL)$, is
also an isomorphism that still satisfies that $\epsilon
^{-1}(t)=\cL_0^{(t)}$ and part (i) of the theorem above.

\end{remark}

\begin{corollary}\label{criterion}
\label{cor}

Let  $A$ be an abelian variety of $\GL_2$-type over $k$ with a
polarization $\cL $ over $k$ of degree $d\geq 1$. Then, $A$ is
principally polarizable over $k$ if and only if there exists $t\in
(R_0)_+$ satisfying $\Norm _{F/\Q }(t)=d$ and $ \cL ^{(t^{-1})}\in
\NS (A_{\bar k})^{G_k}$.
\end{corollary}

{\em Proof. } Assume $\cL ^{(t^{-1})}\in \NS (A_{\bar k})^{G_k}$
for some $t\in (R_0)^*_+$ such that $\Norm _{F/\Q }(t)=d$. Then
$\cL ^{(t^{-1})}$ must be ample because $t$ is totally positive
and principal because $\deg (\cL ^{(t^{-1})})= \deg (\cL ) \Norm
_{F/\Q }(t^{-1}) = 1$. Conversely, if $\cL_0 $ is a principal
polarization on $A$ over $k$, it must be  of the form $\cL _0=\cL
^{(t^{-1})}$ for some $t\in (R_0)_+$ whose  norm from $F$ over
$\Q$ is $d$. $\Box$

The above corollary provides an effective criterion to decide
whether a polarized abelian variety $(A, \cL )$ of real
$\GL_2$-type over $k$ is principally polarizable over $k$. Indeed,
denote by $M$ and $T$ the matrices of the alternating Riemann form
$E$ attached to $\cL$ and $t\in \End_k^+A$ with respect to a fixed
basis of $H_1(A,\Z)$ respectively. It then suffices to check
whether $M\cdot T^{-1}\in \GL_{2g}(\Z)$ for any $t\in (R_0)_+$ (up
to multiplicative elements in $R_0^{* 2}$ if $F=F_0$ or $\Norm
_{F/F_0}(R^*)$ if $[F:F_0]=2$)  such that $\Norm _{F/\Q }(t)=d$ .
Note that if $(A, \cL )$ is an abelian variety of $\GL_2$-type
together with a polarization $\cL$ of primitive type $(1, d_2,
\dots, d_g)$ and degree $d = d_2\cdot ... \cdot d_g$, we only need
to consider those $t\in (R_0)_+ $ of norm $d$ such that $t/m\notin
R_0$ for any $m>1$.

\section{Abelian surfaces of  $\GL_2$-type}\label{sgl2}

\subsection{Principal polarizations on abelian surfaces}
Let $C/k$ be a smooth projective curve defined over a number field
$k$ and let $A = J(C) = \Pic ^0(C)/k$ denote the Jacobian variety
of $C$. The sheaf of sections of the Theta divisor $\Theta _C$
associated to the curve is a principal polarization $\cL (\Theta
_C)$ defined over $k$. For any other curve $C'/k$ such that
$A\stackrel{k}\simeq J(C)\stackrel{k}\simeq J(C')$, Torelli's
theorem asserts that $C$ and $C'$ are isomorphic over $k$ if and
only if there exists $u\in \Aut_{k} A$ such that $u^*(\Theta
_{C'}) = \Theta _C$ in $\NS (A_{\bar k})^{G_k}$.

For an abelian variety $A/k$, this leads us to introduce the set
$T (A_k)$ of $k$-isomorphism classes of smooth algebraic curves
$C/k$ such that $J(C)\stackrel{k}\simeq  A$. Let $\tau (A_k)=|T
(A_k)|$. The Torelli map $C\mapsto \cL (\Theta _C)$ induces an
inclusion of finite sets
$$
\mathrm{T}(A_k)\subseteq \Pi(A_k).
$$

While a principally polarized abelian surface over an algebraic
closed field can only be the jacobian of a curve or the product of
two elliptic curves, the panorama is a little wider from the
arithmetical point of view. We refer the reader to \cite{Mi1} for
an account of Weil's restriction of scalars of abelian varieties
over number fields.

\begin{theorem}
\label{Torelli}

Let $A/k$ be an abelian surface with a principal polarization $\cL
$ defined over $k$. The polarized abelian variety $(A, \cL)$ is of
one of the following three types:

\begin{enumerate}

\item
$(A, \cL ) \stackrel{k}\simeq (J(C), \cL (\Theta _C))$, where
$C/k$ is a smooth curve of genus two.

\item
$(A, \cL )\stackrel{k}\simeq (C_1 \times C_2, \cL _{\mbox{can}})$,
where $C_1$ and $C_2$ are elliptic curves  over $k$ and $\cL
_{\mbox{can}}$ is the natural product polarization on $C_1\times
C_2$.

\item
$(A, \cL )\stackrel{k}\simeq (\Res_{K/k} C, \cL _{\mbox{can}})$,
where $\Res_{K/k} C$ is the Weil restriction of an elliptic curve
$C$ over a quadratic extension $K/k$, and $\cL _{\mbox{can}}$ is
the polarization over $k$ isomorphic over $K$ to the canonical
polarization of $C\times {}^{\sigma}C$.

\end{enumerate}
\end{theorem}

{\em Proof. } It is well known that $(A, \cL )$ is isomorphic over
$\bar k$ to either the canonically polarized Jacobian variety
$(J(C), \cL(\Theta _C))$ of a smooth curve of genus two or to the
canonically polarized product of two elliptic curves
(cf.\,\cite{We2}).

Let us first assume that $(A,\cL)$ is irreducible. We then know
that there exists a curve $C/\qbar $ of genus two such that $(A,
\cL )\stackrel{\qbar}\simeq (J(C), \cL(\Theta _C))$. Since $\cL =
\cL (\Theta _C) \in \NS (A_{\bar k})^{G_k}$, we claim that $C$
admits a $k $-isomorphism onto all its Galois conjugates
${}^{\sigma }C$ for $\sigma \in G_k$. More precisely, if we regard
$C$ as an embedded curve in $\Pic ^1(C)$, then ${}^{\sigma }C=C +
a_{\sigma }$ for some $a_{\sigma }\in \Pic ^0(C)(\qbar )$. Indeed,
this follows from the fact that the sheaves ${}^{\sigma}\cL=\cL
(\Theta _{\,{}^{\sigma }C})$ are all algebraically equivalent and
$h^0(\cL )=1$. In particular, we obtain that the field of moduli
$k_C$ of $C$ is contained in $k$. Let us now show that $C$ does
admit a projective model over $k$. We distinguish two cases
depending on whether the group $\Aut C$ is trivial or not.

If $\Aut C\not \simeq \Z /2 \Z $, Cardona and Quer have recently
proved that $C$ always admits a model over its field of moduli
(cf. \cite{CaQu}).

We now consider the case that the hyperelliptic involution $v$ on
$C$ generates the group of the automorphisms of $C$. Then, as
shown by Mestre in \cite{Me}, there is a projective model $C/K$ of
$C$ over a quadratic extension $K/k$. Let $\sigma \in G_k$ such
that $\sigma$ does not act trivially on $K$. There is an
isomorphism $\varphi _{\sigma }:C\stackrel {\simeq }{\ra
}{}^{\sigma }C$ of $C/ {K}$ onto $C^{\sigma }$ given by the
translation by $a_{\sigma }$ map on $\Pic ^1(C)$. The map $\Pic
^1(C)\ra \Pic ^1(C)$, $D\mapsto D + a_{\sigma }+
{}^{\sigma}a_{\sigma }$ descends to an automorphism of $C$,
${}^{\sigma }\varphi _{\sigma }\circ  \varphi _{\sigma }$, which
cannot be the hyperelliptic involution $v$, since $v=-1_{J(C)}$ on
$J(C)$. As  we are assuming that $\Aut (C)\simeq \Z / 2 \Z $, we
obtain that ${}^{\sigma }\varphi _{\sigma }= \varphi _{\sigma
}^{-1}$. Now, Weil's criterion on the field of definition of an
algebraic variety applies to ensure that $C$ admits a projective
model over $k$ (cf. \cite{We}).

Now, assume that $(A,\cL)$ is reducible over $\overline k$, i.e.,
$A\stackrel{\overline k}\simeq C_1\times C_2$ for some elliptic
curves $C_1$ and $C_2$. If both $C_1, C_2$ and the isomorphism are
defined over $k$, then $(A, \cL )\stackrel{k}\simeq (C_1 \times
C_2, \cL _{\mbox{can}})$. Otherwise, $C_i/K$ must be defined over
a quadratic extension $K/k$ and  $C_1={}^{\sigma }C_2$ where $\Gal
(K/k)= \langle \sigma \rangle $, since the product $C_1\times C_2$
is defined over $k$. This is equivalent to saying that
$A\stackrel{k}\simeq \Res _{K/k}C_1$. $\Box$

\begin{definition}
We say that a principal polarization $\cL $ on an abelian surface
$A$ over $k$ is {\em split} if $(A, \cL )\stackrel {\bar k
}{\simeq }(C_1\times C_2, \cL _{\mbox{can}})$ for some elliptic
curves $C_i/\bar k$.  We shall denote  by  $\sigma (A_k)$ the
number of $k$-isomorphism classes of split principal polarizations
on $A$ over $k$.
\end{definition}

\begin{corollary}

Let $A/k$ be an abelian surface. Then

$$
\pi (A_k) = \tau (A_k) + \sigma (A_k).
$$

\end{corollary}

It may very well be the case that $\sigma (A_k)\geq 2$ for some
abelian surface $A/k$. This amounts to saying that $A\simeq
C_1\times C_2\simeq C_3\times C_4$ as unpolarized abelian
varieties for two different nonordered pairs of elliptic curves
$(C_1, C_2)$ and $(C_3, C_4)$ (cf. \cite{BiLa}, p. 318).

\subsection{Principal polarizations on abelian surfaces of  $\GL_2$-type}

We now focus our attention on abelian surfaces of  $\GL _2$-type.
As an immediate consequence of Theorem \ref{pi}, we obtain the
following.

\begin{corollary}\label{surpol}

Let $A$ be a principally polarizable
abelian surface over $k$ such that $\End _k A=R$ is an order in
a quadratic field $F$. Then

$$
\pi(A_{k})=\begin{cases}
    2 & \text{if $F$ is real and all units in $R$ have positive norm,}\\
    1 & \text{otherwise.}

  \end{cases}
$$

\end{corollary}

Since  every abelian variety of $\GL_2$-type over $k$ is simple
over $k$, note that the second possibility of Theorem
\ref{Torelli} cannot occur when the surface is of this type. The
following definition was first introduced in \cite{El}.

\begin{definition}
An elliptic curve $C/\bar k$ is a $k$-curve if $C$ is isogenous to
all its Galois conjugates $C^{\sigma }$, $\sigma \in \Gal (\bar
k/k)$. A $k$-curve $C$ is completely defined over an extension
$K/k$ if $C$ is defined over $K$ and it is isogenous over $K$ to all
its Galois conjugates.
\end{definition}

\begin{lemma}\label{k-curve}

Let $A$ be an abelian surface of real $GL _2$-type over $k$, which
is $k$-isogenous to the Weil restriction $\Res_{K/k} C$ of an
elliptic curve $C/K$, where $K/k$ is a quadratic extension. Then
$C$ is a $k$-curve completely defined over $K$, its $j$-invariant
$j(C)\notin k$ and $\Q\otimes\End_K A$ is isomorphic to
$\M_2(\Q)$.
\end{lemma}

{\em Proof. }  Let $\sigma$ be the nontrivial automorphism of
$\Gal (K/k)$. We have that $\Q \otimes \End_K C=\Q$ and there is
an isogeny $ \mu: C\rightarrow {}^{\sigma} C$ defined over $K$,
since otherwise $F=\Q \otimes \End_k A$ must be either $\Q$ or
contain an imaginary quadratic field. Therefore, $C$ is a
$k$-curve completely defined over $K$ and $\mu \circ
{}^{\sigma}\mu\in\Q$.  Then, $\Q\otimes\End_K A\simeq \M_2(\Q)$
and $ \mu \circ {}^{\sigma}\mu=\pm \deg \mu$. Since $F$ is a
totally real quadratic field, $\mu\circ {}^{\sigma}\mu=\deg \mu$
and $F=\Q(\sqrt{ \deg \mu})$. Consequently, $\mu$ cannot be an
isomorphism and hence $j(C)\notin k$. $\Box$

\begin{proposition}\label{resweil}

Let $C_1$, $C_2$  be two elliptic curves defined over the
quadratic fields $K_1$, $K_2$ respectively. If $\Res_{K_1/\Q}C_1$
and $\Res_{K_2/\Q}C_2$ are $\Q$-isomorphic and they are of real
$\GL_2$-type over $\Q$, then $K_1=K_2$ and $C_1$ is isomorphic
over $K_1$ either to $C_2$ or its Galois conjugate.

\end{proposition}

{\em Proof. } Let $A=\Res_{K_i/\Q} C_i$. First, we show that
$K_1=K_2$. If $C_i$ does not have complex multiplication, we then
know by the previous lemma that $\Q\otimes
\End_{K_i}A=\Q\otimes\End _{\overline{\Q}}A\simeq M_2(\Q)$ and
hence $\Q\otimes\End _{K_1\cap K_2}A\simeq M_2(\Q)$. It follows
that $K_1=K_2$.

Assume then that both $C_1$ and $C_2$ have complex multiplication
by the same  imaginary quadratic field $L$ because $C_1$ and $C_2$
must be  isogenous. Again by the lemma, $K_i=\Q (j(C_i))$. We have
that $j(C_i)= j(\ga_i )$ for  an invertible ideal $\ga_i$ of the
order $\End (C_i)$ of $L$. Since $K_i$ is quadratic,  the ideal
classes of  $\ga_i$ and $\ga_i^{-1}$ are equal and, therefore,
$j(\ga_i)=j(\ga_i^{-1})=j(\overline{\ga_i})=\overline{j(\ga_i)}$,
where $\overline{\phantom{c}}$ denotes the complex conjugation. It
follows that   $K_i=\Q (j(C_i))$ is a real quadratic field. Let us
assume that $K_1\neq K_2$ and let $\Gal (K_i/\Q )=\langle \sigma
_i\rangle $. Since $C_1\times C_1^{\sigma _1}\simeq C_2\times
C_2^{\sigma _2}$ over $K_1\cdot K_2$ but not over $K_i$, there
exists an isogeny $\mu: C_1\rightarrow C_2$ defined over $K_1\cdot
K_2$ but not over $K_i$. Let $\tau \in \Gal (K_1\cdot K_2/\Q)$ be
the automorphism which does not act trivially over $K_1$ or over
$K_2$. Then, $\widehat{\mu}=\mu\times {}^{\tau}\mu\in\Q\otimes
\End _{K_1\cdot K_2}A\setminus \Q\otimes \End_{K_i}A$. Let us
denote by $\cA$ the $\Q$-algebra generated by $\Q\otimes
\End_{K_i}A$ and $\widehat{\mu}$. From the inclusions
$$
\Q\otimes \End_{K_i}A\simeq M_2(\Q)\subsetneq
\cA \subseteq \Q\otimes\End_{\overline{\Q}}A\simeq M_2(L)\,,
$$
we obtain that $\cA = \Q\otimes \End_{\overline{\Q}}A$ and, thus,
all endomorphisms of $A$ are defined over $K_1\cdot K_2$. In
particular, we obtain that $\Q\otimes\End_{K_1\cdot
K_2}(C_i)\simeq L$, but this leads to a contradiction since the
totally real field $K_1\cdot K_2$ cannot contain $L$. Therefore,
$K_1=K_2$.

Set $K=K_i$, $C=C_1$ and $\Gal (K/\Q)=\langle \sigma\rangle$. Now,
we will prove that $C$ is isomorphic either to $C_2$ or to
${}^{\sigma}C_2$  over $K$. Let us denote by $\pi :A\rightarrow C$
the natural projection over $K$,  so that $\ker \pi={}^{\sigma}C$.
Then the period lattice of $C$  is  $\Lambda=\{ \int_{\gamma}
\pi^{*} (\omega) \mid \gamma \in H_1(A,\Z) \}$, where $\omega$ is
the invariant differential of $C$.  Given $\alpha\in \Q\otimes
\End_K A$, we will denote by $\Lambda_{\alpha}$ the set $\{
\int_{\gamma} \alpha^{*}(\pi^{*} (\omega) )\mid \gamma \in
H_1(A,\Z) \}$. When $\alpha^{*}(\pi^{*} (\omega))\neq 0$,
$\Lambda_{\alpha}$ is the period lattice of  a certain  elliptic
curve over $K$; these lattices cover all the $K$-isomorphism
classes of  elliptic curves which are optimal quotients of $A$
over $K$.  We will prove  that all these classes are also obtained
when $\alpha$ only runs over $F^*$. Let $w\in\End_KA$ be the
composition of the morphisms
$$
A\stackrel{\pi}\rightarrow C
\stackrel{i}\hookrightarrow A=C\times {}^{\sigma}C\,,
$$
where $i$ is the natural inclusion. We have  $w^{*}
H^0(A,\Omega^1_{A/K})= K \pi^{*} (\omega)$ and  $w^2=w$. Moreover,
$\Q\otimes \End_KA=F\oplus w\cdot F$ since $w\notin F$ and $
\Q\otimes \End_K A$ is a $F$-algebra of dimension $2$. Now, it
suffices  to use the equality  $w\cdot \Q\otimes \End_KA= w\cdot
F$. Since  for every integer $m\neq 0$, the classes  corresponding
to $\Lambda_{\alpha}$ and $\Lambda_{m\alpha}$  are isomorphic, we
can assume that  $\alpha\in \End_{\Q} A$ and in this case
$\Lambda_{\alpha}$ is a sublattice of $\Lambda$.

There exists a cyclic isogeny between $C$ and ${}^{\sigma}C$ over
$K$ of a certain degree $n$, which extends to an endomorphism
$\beta \in \End_{\Q}A$ with the following properties:

\begin{enumerate}
\item
$\beta$ restricted to ${}^{\sigma} C$ (resp. $C$) provides a
cyclic isogeny between ${}^{\sigma} C$ and $C$ (resp. $C$ and
${}^{\sigma} C$) of degree $n$.

\item $\Lambda_{\beta}$ is the
period lattice of the elliptic curve  isomorphic to ${}^{\sigma}C$
over $K$ which satisfies $\Lambda/\Lambda_{\beta}\simeq \Z/n\Z$.
\item  $\beta^2=n$, $F=\Q( \beta)$ and, moreover,
for all integers $m>1$ we have that $\beta/m\notin \End_{\Q}A$.

\item Due to the previous step, $\beta$ also provides a cyclic
isogeny between ${}^{\sigma} C_2$ and $C_2$ of degree $n$.

\item For all integers $a$, $b$ we have
$\Lambda_{a+b\beta}=a\Lambda+ b\Lambda_{\beta}$.
\end{enumerate}
Let $\Lambda_{\alpha}$ be a lattice corresponding to the
$K$-isomorphism class of $C_2$ and write $\alpha=a+ b\beta$ with
$a$, $b$ integers which we can assume to be coprime. We can also
assume that $a,b\neq 0$, since otherwise the statement is obvious.
Set $d=\gcd (a, n)$.  Using the fact
$\Lambda/\Lambda_{\beta}\simeq \Z/n\Z$, it can be checked that
$\Lambda_{\alpha}=a\Lambda+ b \Lambda_{\beta}$ is the  sublattice
$d \Lambda+  \Lambda_{\beta}$ of $\Lambda$ of degree $d$. Then, we
can take  $a=d$ and $b=1$. Now, we  can see that
$\Lambda_{\alpha\cdot \beta}=n \Lambda+ d\Lambda_{\beta}$ is a
sublattice of $\Lambda_{\alpha}$ of degree $d\cdot n$. Therefore,
$d$ must be $1$ and $C$ and $C_2$ are isomorphic over $K$. $\Box $

\begin{corollary}\label{number} Let $A$ be a principally polarizable
abelian surface of $GL _2$-type over $\Q $. Then, $\sigma
(A_{\Q })\leq 1$ and, in particular,
$$
\tau (A_{\Q})=\begin{cases}
    \pi(A_{\Q})-1 & \text{if $A \stackrel{\Q}\simeq \Res_{K/\Q}C$ for a
    $\Q$-curve $C/K$,
    } \\
    \pi(A_{\Q})& \text{otherwise.}
  \end{cases}
$$
\end{corollary}

\begin{remark}

Note that if $\End _{\bar \Q}A=\End _{ \Q}A$, then
$\pi(A_{\Q})=\pi(A_{\bar \Q})$ and $\tau (A_{\Q})= \tau (A_{\bar
\Q})$. Furthermore, if $\tau (A_{\Q})=2$ then the two curves over
$\Q$ which share $A$ as Jacobian are nonisomorphic over $\bar \Q$.

\end{remark}

Now, Theorem \ref{main} is a consequence of corollaries
\ref{surpol} and \ref{number}.

\begin{corollary}

There are infinitely many genus two curves over $\Q$ such that
their Jacobians are of real $\GL_2$-type and isomorphic to the
product of two elliptic curves as unpolarized abelian varieties.

\end{corollary}

{\em Proof. } As a consequence of Proposition \ref{surpol} and Theorem
\ref{number}, we have  that for every quadratic $\Q $-curve $C$ such that
$A=\Res_{K/\Q} C$ is of real $\GL_2$-type and that the order $\End_{\Q } A$
contains no units of negative norm, there is a genus two curve over $\Q$ whose Jacobian is
isomorphic to $A$ as unpolarized abelian varieties.
Fix $d = 3$ or $7$. It is known
that for every quadratic field $K$ there exists a  $\Q $-curve $C'/K$ without CM
such that $K = \Q( j(C'))$ and an isogeny from $C'$ onto its
Galois conjugate of degree $d$. For
such a curve there is an isomorphic curve $C/K$ with $A=\Res_{K/\Q} C$ being of
real $\GL_2$-type if and only $d$ is a norm of $K$ (see \cite{Qu}). Now, the statement follows
from the fact that $\Z [\sqrt d]$ contains no units of negative norm and the
existence of infinitely many quadratic fields $K$ such that
$d\in \Norm _{K/\Q }(K^*)$.
$\Box$

In subsection \ref{constructing}, we will show how for some
quadratic $\Q$-curves Theorem \ref{main} allows us to construct a
genus two curve over $\Q$ whose Jacobian is isomorphic over $\Q$
to its Weil's restriction.

In view of the above, it would be useful to have a criterion for
deciding whether an abelian surface of real $\GL _2$-type is
principally polarizable or not. In addition to Corollary
\ref{criterion}, we now characterize the existence of a principal
polarization on an abelian surface of real $\GL_2$-type under some
arithmetical restrictions on the ring of the endomorphisms.

\begin{proposition}\label{principal}
Let $A$ be an abelian surface of $\GL _2$-type over $k$ by a real
quadratic field $F$ of class number $h(F)=1$ and assume that $\End
_kA=\cO $ is the ring of integers of $F$. Then

\begin{enumerate}
\item
If $\cO^*$ contains some unit of negative norm, then $A$ is
principally polarized over $k$ and, in particular, the degree of
any polarization on $A$ over $k$ is a norm of $F$.

\item
If $\cO^*$ contains no units of negative norm and $\cL $ is a
polarization of degree $d$ over $k$, then either $d$ or $-d$ is a
norm of F. In the first case, $A$ is principally polarizable over
$k$.
\end{enumerate}
\end{proposition}

{\em Proof. } First, we will prove that there is an invertible
sheaf on $A$ in $\NS (A_{\bar k})^{G_k}$ of degree $1$. Since $A$
is defined over $k$, there is a polarization $\cL$ on $A$ defined
over $k$. Let us denote by $E$ the corresponding alternating
Riemann form on $H_1(A,\Z)$. Since $H_1(A,\Z) $ is a free
$\cO$-module of rank $2$ and $h(F)=1$, we have that $H_1(A,\Z) =
\cO \gamma_1 \oplus  \cO \gamma_2$ for some $\gamma_1,\gamma_2\in
H_1(A,\Z) $. Consider the morphism of groups
$$
\phi:\cO \longrightarrow \Z\,, \quad \alpha \mapsto E(\alpha \gamma_1,\gamma_2)\,.
$$
By the nondegeneracy of the trace, there is $\delta \in F$ such
that $\phi (\alpha)=\Tr_{F/\Q} (\alpha \cdot \delta)$. Due to the
fact that $\delta$ lies in the codifferent of $\cO$,  there is
$v\in\cO$ such that $\delta= v/\sqrt {\Delta}$, where $\Delta$ is
the discriminant of $F$. Set $\cL_0=\cL^{(v^{-1})}\in\Q\otimes
\NS(A_{\bar k})^{G_k}$ and denote by $E_0$ its alternating Riemann
form. From the relations
$$
E_0(\alpha \gamma_1,\gamma_2)= \Tr_{F/\Q}(\alpha/\sqrt{\Delta})\in\Z\,,\quad
E_0(\gamma, \alpha \gamma)=0
$$
for all $\alpha \in \cO$ and $\gamma \in H_1(A,\Z)$, we obtain
$\cL_0\in \NS(A_{\bar k})^{G_k}$. Moreover, taking a basis
$\gamma_1, \alpha \gamma_1,\gamma_2,\alpha \gamma_2$ of
$H_1(A,\Z)$ for a suitable $\alpha$, an easy check shows that the
matrix of $E_0$ with respect to this basis has determinant $ 1$.

Now, assume  that $\cO^*$ contains some unit $u$ of negative norm.
Then, either $\cL_0$, $\cL_0^{(-1)}$, $\cL_0^{(u)}$ or
$\cL_0^{(-u)}$ must be a principal polarization. By Corollary
\ref{criterion}, we also obtain that the degree of any
polarization on $A$ over $k$ is a norm of $F$.

In the case $-1\not \in \Norm _{F/\Q }(\cO ^*)$,  by remark
\ref{degree1}, there is $t\in \cO$ such that $\cL_0^{(t)}=\cL$,
with $\Norm_{F/\Q}(t)=\pm \deg \cL$. If $\deg \cL$ is a norm of
$F$, then $\Norm_{F/\Q}(t)=\deg \cL$. Therefore, either $t$ or
$-t$ lies  in $\cO_+$ and, thus, either $\cL_0$ or $\cL_0^{(-1)}$
is a polarization. $\Box$

\section{Explicit examples}\label{examples}

In the previous section we have shown the different possibilities
for abelian surfaces of $\GL_2$-type over $\Q$. We now illustrate
these possibilities  with explicit examples for the case of real
$\GL_2$-type. We will not consider abelian surfaces of CM $\GL
_2$-type because they are less interesting. Indeed, these abelian
surfaces have a unique primitive polarization $\cL$ and they are
principally polarizable if and only if $\cL $ is principal.
Assuming Serre's conjecture to be true, we must look for these
examples among the abelian varieties $A_f$ attached by Shimura to
an  eigenform $f\in S_2(\Gamma_0(N))$. In order to present these
abelian varieties as Jacobians of curves, we will use the
procedure described in \cite{GoGoGu}, which is based on jacobian
Thetanullwerte (see \cite {Gu}).

All the computations were performed with {\tt Magma v.2.7} using
the package {\tt MAV} written by E. Gonz\'alez and J. Gu\`ardia
(\cite{GoGu}). Both the package and some files required to
reproduce them are available via the web pages of the authors.

\subsection{Modular Case}\label{mod} We summarize here
some facts about modular abelian varieties, fixing also the
notation used in the examples in the following subsections. Let
$f=\sum_{n\geq 1} a_n q^n$ be a normalized newform of
$S_2(\Gamma_0(N))$ where, as usual, $q=e^{2\pi iz}$. Attached to
$f$, Shimura constructed an abelian variety $A_f$ over $\Q$ in two
different ways. In  \cite{Shi71}, it is presented this abelian
variety as a subvariety of $J_0(N)$ while in \cite{Shi73} it is
constructed as an optimal quotient of $J_0(N)$ and it is proved
that both are dual. We will take $A_f$ as subvariety, and this
will not be a restriction because we only will consider the case
principally polarized and in this situation both abelian varieties
are isomorphic (over $\Q)$. More precisely, let $\mathbf T$ be the
Hecke algebra of endomorphisms of the Jacobian variety $J_0(N)$ of
$X_0(N)$, and let $I_f$ be the kernel of the map $ {\mathbf T} \to
\Z[a_1,a_2,\dots]$ which identifies every Hecke operator with the
corresponding eigenvalue of $f$. Then $A=J_0(N)/I_f J_0(N)$ is the
abelian variety attached by Shimura as an optimal quotient. We
recall that $K_f=\Q(\{a_n\})$ is a number field of degree
$n=[K_f:\Q]=\dim A$,  the endomorphism algebra $\Q\otimes
\End_{\Q}A$ is the $\Q$-algebra generated by $\mathbf T$ acting on
$A$, $\Q\otimes \mathbf T/I_f$,   and this is isomorphic to $K_f$.
Let us denote by $\pi:J_0(N) \rightarrow A$ the natural projection
over $\Q $. There is a $\Z$-submodule $H$ of $H_1(J_0(N),\Z)$ of
rank $2 n$ such that $H_1( J_0(N),\Z)= \ker \pi_{*} \oplus H$.
Note that $\ker \pi_{*}=I_f H_1( J_0(N),\Z)$.

It is well known that $\pi^{*} H^0(A,\Omega^1_{A/\C})$  is the $\C$-vector space generated by
the Galois conjugates of $f(q)\, dq/q$ and $\pi^{*}
H^0(A,\Omega^1_{A/\Q})$ is the subspace obtained by taking the
modular forms with rational $q$-expansion.
For a fixed rational basis $h_1,\dots, h_n$ of $\pi
^{*}H^0(A,\Omega^1_{A_f/\Q})$, the Abel-Jacobi map induces an
isomorphism of complex torus:
$$
A (\C ) \longrightarrow \C^n/\Lambda\,,\quad  P\mapsto (\int_0^P
h_1,\dots \int_0^P h_n)\,,
$$
where
$$\Lambda= \{ (\int_{\gamma} h_1,\dots \int_{\gamma} h_n)| \gamma\in H_1(X_0(N),\Z)\}=\{ (\int_{\gamma} h_1,\dots \int_{\gamma} h_n)| \gamma\in H\}\,.$$
The abelian variety $A_f$, viewed as subvariety of $J_0(N)$, is described by
$$
\begin{array}{rcl}
A_f (\C )&\longrightarrow &\C^n/\Lambda_f\\
 P&\mapsto & (\int_0^P h_1,\dots
\int_0^P h_n)\,,
\end{array}
$$
where $\Lambda_f= \{ (\int_{\gamma} h_1,\dots \int_{\gamma} h_n)| \gamma\in H_f\}$ and $H_f=\{\gamma\in H_1( X_0(N), \Z)\mid I_f \gamma=\{ 0\}\}$. Note $\Lambda_f$ is a sublattice of $\Lambda$.
Obviously, given $T\in\Q \otimes \mathbf T$, we have that $T\in \End_{\Q} A_f$
if and only if $T$ leaves $\Lambda$ stable or equivalently $T$ leaves $H_f$ stable. Let $\Theta$ be the canonical
polarization on $J_0(N)$ and $E$ its corresponding Riemann form, which is
obtained from the intersection numbers on $H_1(X_0(N),\Z)$. From now on, we shall call
the polarization $\cL $ on $A_f$ obtained from the canonical polarization $\Theta$ on $J_0(N)$ the {\em canonical polarization} on $A_f$. Note
that the corresponding Riemann form  $E_{\Lambda}$ on $\Lambda $
is obtained as the restriction of $E$ to $\Lambda$ and that, although $\Theta$ is
principal, $\cL$ may not be.

We know by Ogg (see \cite{Og}) that the cusps of $X_0(N)$
associated to $1/d$ for $d\mid N$ and $\gcd (d,N/d)=1$ are
rational points on $X_0(N)$. Moreover, the divisors of degree $0$
generated by these cusps are torsion points on $J_0(N)(\Q)$ and
the same holds for their projections on $A_f$, although their
orders can decrease. For a given torsion  $\Gal
(\overline{\Q}/\Q)$-stable subgroup $G$ of $A_f$, the abelian
variety $A_f/G$ is isogenous to $A_f$ over $\Q $. Moreover, both
abelian varieties can be isomorphic over $\Q$ only if $G$ is the
kernel of some endomorphism in $\End_{\Q}(A_f)$. As complex tori,
if we identify $A_f(\C )=\C ^g/\Lambda_f $, then $A_f/G(\C
)=\C^g/\Lambda_G$, where $\Lambda_G=\langle \Lambda_f , G\rangle
$. The polarization $\cL$  induces a polarization $\cL_G$ on
$A_f/G$, whose alternating Riemann form $E_G$ is given by $E_G=\#
G\cdot E_{\Lambda}: \Lambda _G\times \Lambda _G\rightarrow \Z $
and its degree is $\# G\cdot \deg \cL$. Notice that
$\End_{\Q}(A_f)$ and $\End_{\Q} (A_f/G)$ can be different.

\subsection{Nonisomorphic genus two curves with $\Q$-isomorphic Jacobian}
We begin by studying the unique two-dimensional factor $S_{65,B}$ of
$J_0(65)$, given by the newform
$$
f= q + aq^2 + (1-a)q^3 +q^4 -q^5  +(a-3)q^6\dots
$$
with $a=\sqrt{3}$. We know that $\End_{\Q}(A_f)=\Z [\sqrt{3}]$ and
moreover $A_f$ is simple in its isogeny class over
$\overline{\Q}$, because $f$ does not admit any extra-twist.
Hence, each principal polarization on $A_f$ is the sheaf of
sections of the Theta divisor of a smooth curve $C$ of genus two
such that $A_f\stackrel{\Q}\simeq J(C)$. A basis of  the
$\Z$-module $H_f$ is spanned by the  modular symbols
$$
\begin{array}{l}
 \gamma_1=\{-\frac{1}{15}, 0\}  -\{-\frac{1}{30},
0\}+\{-\frac{1}{40},0\}-\{-\frac{1}{60},0\}, \\
\gamma_2=    \{-\frac{1}{20}, 0\}  -\{-\frac{1}{35},
0\}+\{-\frac{1}{50},0\}-\{-\frac{1}{55},0\},\\
\gamma_3=    \{-\frac{1}{15}, 0\} -\{-\frac{1}{26}, 0\} +
\{-\frac{1}{40},0\}-\{-\frac{1}{50},0\}-\{-\frac{2}{5},-\frac{5}{13}\}, \\
\gamma_4=    \{-\frac{1}{30}, 0\}  -\{-\frac{1}{45}, 0\} + \{-\frac{1}{52}, 0\}
-\{-\frac{1}{55}, 0\}+\{-\frac{2}{5},-\frac{5}{13}\}.
\end{array}
$$
An integral basis  of $H^0(A_f,\Omega^1)$ is given by the forms $h_1=({}^{\sigma}f-f)/\sqrt{3}$, $h_2=(f+
{}^{\sigma}f)/2$. By
integrating these differentials  along the paths $\gamma_1, \dots
,\gamma_4$, we obtain an analytic presentation of the abelian
surface $A_f$ as a complex torus $\C^2/\Lambda $.  The
restriction of $E$ to $H_f$ is the Riemann form of a polarization
on $A_f$, given by the following Riemann matrix:
$$
M_E=\left( E(\gamma_i,\gamma_j)\right)_{i,j}=
\left(    \begin{matrix}
  0 &  2 & -2 & 0\\
 -2 &  0 & -2 & 2\\
  2 &  2 &  0 & 0\\
  0 & -2 &  0 & 0\\
\end{matrix}\right).
 $$
The type of this polarization is $(2,2)$. Thus, the primitive
polarization $\cL_0 $ associated to it is principal. In
\cite{GoGoGu} it has been checked that $(S_{65, B}, \cL_0 )$ is the
polarized Jacobian of the hyperelliptic curve
$$
C_{65,B,1}: Y^2\,=-X^6 - 4X^5 + 3X^4 + 28X^3 - 7X^2 - 62X + 42,
$$
whose absolute Igusa invariants are
$$
\displaystyle \{i_1,i_2,i_3\}=\{ -\frac{2^{6}5\cdot 313^5}{ 13^3},
\frac{139\cdot313^3 701  }{5\cdot 13^3 }, \frac{7\cdot 313^2\cdot59104229}
{2^3 5^2 13^3}\}.
$$
Following Theorem \ref{pi}, we now build a second
polarization on $S_{65,B}$, which will exhibit $S_{65,B}$ as the
Jacobian of a second curve $C_{65,B,2}$ nonisomorphic to
$C_{65,B,1}$ over $\qbar $. Let us consider the Hecke operator
$u=2+T_2=2+\sqrt{3}\in \End_{\Q } A_f$.
This is a non-square totally
positive unit in the ring of integers $\cO$ of $K_f$. The action
of $u$ on $H_1(A_f,\Z)$ with respect to the basis
$\gamma_1,\gamma_2, \gamma_3,\gamma_4$ is given by:
$$
M_u=
\left(    \begin{matrix}
  3 &  2 & -1 & 1\\
  0 &  0 &  1 & 0\\
  0 & -1 &  4 & 0\\
  2 &  1 &  1 & 1\\
\end{matrix}\right).
$$
We now take the polarization $\cL_0^{(u)}$, which is given by  the
Riemann form $E_u$ whose alternating Riemann matrix is $M_{E_u}= 1/2\cdot M_E\cdot M_u$. It is
 a principal polarization on $A_f$ and  nonisomorphic to $\cL_0 $. Since $A_f$ is simple (over $\bar\Q$),   the polarized abelian variety $(A_f,\cL _0^{(u)})$ is the Jacobian of a curve $C_{65,B,2}$ defined over $\Q$.
A symplectic basis with respect to $E_u$ is
$$
\left(
\begin{array}{c}
\delta_1 \\
\delta_2 \\
\delta_3 \\
\delta_4
\end{array}
\right)=\left(
\begin{array}{rrrr}
 1 & 0 & 0 & 0\\
 0 & 3 & 1 & 0\\
 0 & 1 & 0 &-2\\
 0 & 0 & 0 & 1
\end{array}\right)
\left(
\begin{array}{c}
\gamma_1 \\
\gamma_2 \\
\gamma_3 \\
\gamma_4
\end{array}
\right).
$$
In order to identify the polarized abelian variety
$(A_f,\cL^{(u)})$, we compute the period matrix
$\Omega=(\Omega_1\, | \Omega_2)= ((\int_{\delta_k} h_jdq/q)_{k=1,2}
|(\int_{\delta_k} h_j dq/q)_{k=3,4}) $
 of the forms
$h_1, h_2$ along these paths.
Let $Z=\Omega_1^{-1}\Omega_2\in {\mathcal{H}}_2$, which belongs to
the Siegel upper half space. We can now
apply the method of \cite{GoGoGu} to recover an equation of this
curve. We obtain
$$
\begin{array}{rl}
C_{65,B,2}:\, Y^2=&-(X^2+3X+1) (7X^4 +37X^3+ 71X^2 + 44X + 8).
\end{array}
$$
The absolute Igusa invariants
of this curve are
$$
\displaystyle \{i_1,i_2,i_3\}=\{
-\frac{2^{3}3109^5}{5^3 13^4},
\frac{3109^{3} 206639  }{2\cdot 5^3 13^4},
\frac{7\cdot 3109^2 123916753}{2^3 5^3 13^4}\}.
$$
and this makes it clear that the curves $C_{65,B,1}$ and
$C_{65,B,2}$ are nonisomorphic over $\qbar $, although their
Jacobians are isomorphic as unpolarized abelian varieties over $\Q
$.

\vskip 0.5cm
\noindent
\subsection{A  genus two curve with Jacobian isomorphic to the Weil
restriction of a quadratic $\Q$-curve}
Let us now consider the factor $A_f$ of $J_0(63)$ corresponding to the
newform
$$
f=q+aq^2+q^4-2aq^5+q^7-aq^8+6q^{10}+2q^{11}+2q^{13}+\dots ,
$$
where $a=\sqrt{3}$. Again  $\End_{\Q}(A_f)\simeq \Z [\sqrt
{3} ]=:\cO$ is the integer ring of $K_f=\Q(\sqrt3)$. Let $\cL_0 $ be the primitive canonical polarization on $A_f$
induced from $J_0(63)$. It turns out that $\cL_0 $ is principal and
that the polarized abelian variety $(A_f,\cL )$ is the Jacobian of
the hyperelliptic curve (see \cite{GoGoGu}):
$$
C_{63,B}:\quad Y^2=\, -3 X^6 + 162 X^3 + 81.
$$

We have that $f$ does not have complex multiplication and  the quadratic character of
 $L=\Q (\sqrt {-3})$ is an
extra-twist for $f$.
 Since the discriminant of
$L$ is a norm of $K_f$, we know (see \cite{GoLa}) that, up to Galois conjugation,
there is a unique quadratic $\Q $-curve $C$ defined over $L$ such
that $C$ is an optimal quotient of $A_f$ defined over $L$.
In fact, the  invariant differentials of the optimal quotients of $A_f$, when pulled back to $H^0(A_f, \Omega^1_{A_f/L})$ are given by
$
\langle \left (\frac{1-i}{2} (a+ b\sqrt 3) f+  \frac{1+i}{2} (a-
b\sqrt 3) {}^{\sigma}f \right ) dq/q \rangle \,,
$
with $a+ b \sqrt 3$ running over $ K_f^{*}$.

In order to obtain the other principal polarization on  $A_f$, we  proceed as in the previous example. Let $u=T_2+T_{13}=\sqrt{3}+2\in \cO ^*\setminus \cO ^{\ast 2}$ be the \newline
Hecke operator acting on $A_f$ and let $\cL_0^{(u)}$ be the principal
polarization on $A_f$ associated to it by Proposition \ref{pol}.
A symplectic basis of $H_1(A_f,\Z)$ with respect the Riemann form $E_u$ associated to $\cL_0^{(u)}$
is:
$$
\begin{array}{rl}
_{\gamma_1=}& _{
  \{-\frac1{24}, 0\} -\{-\frac1{28}, 0\} + \{-\frac1{30}, 0\}
  -\{-\frac1{51}, 0\} -\{-\frac13, -\frac27\},  }\\
_{\gamma_2=}&    _{  2(\{-\frac1{24}, 0\}  -\{-\frac1{28}, 0\}+ \{-\frac1{39}, 0\}
  -\{-\frac16, -\frac17\}) -5\{-\frac1{57}, 0\} }\\
&  _{+3(-\{-\frac1{36}, 0\}
  + \{-\frac1{49}, 0\} -\{-\frac1{51}, 0\}
  + \{-\frac1{54}, 0\}
  +\{-\frac1{60}, 0\} + \{-\frac1{3}, -\frac27\}),}   \\
 _{\gamma_3=}& _{  -\{-\frac1{28}, 0\} + \{-\frac1{36}, 0\} + \{-\frac1{45}, 0\}
 -\{-\frac1{49}, 0\} + \{-\frac1{51}, 0\}
 -\{-\frac1{54}, 0\}
  -\{-\frac16, -\frac17\} + \{\frac37, \frac49\},    }\\
  _{\gamma_4=}& _{  2(-\{-\frac1{36}, 0\}
 + \{-\frac1{49}, 0\}-\{-\frac1{51}, 0\} + \{-\frac1{54}, 0\}
 -\{-\frac1{57}, 0\})
  +\{-\frac1{24}, 0\} }
\\
& _{
 +  \{-\frac1{39}, 0\}
-\{-\frac1{45}, 0\}
   +
    \{-\frac1{60}, 0\} + \{-\frac13, -\frac27\}  -\{\frac37, \frac49\}.}
\end{array}
$$
We take the basis of $H^0(A_f,\Omega ^1_{A_f/L})$ given by
$$
\begin{array}{rl}
g_1=&\left ( \frac{1-i}{2} (1+ \sqrt 3) f+
\frac{1+i}{2} (1- \sqrt 3) {}^{\sigma}f \right ) dq/q
\\ \\
=&
\left(  (1-i\sqrt3)q +(3-i\sqrt3)q^2+(1-\sqrt3)q^4+\dots\right)dq/q
\end{array}
$$
and its conjugate $\overline{g_1}$.
Let $\Omega=(\Omega_1,\mid \Omega_2)$ the period matrix of this basis $g_1, \overline{g_1}$
with respect to $\gamma_1, \gamma_2, \gamma_3, \gamma_4$ and
take $Z=\Omega_1^{-1}\Omega_2$.
When we compute the even Thetanullwerte
corresponding to   $Z$,  we find that, up to high
accuracy, exactly one of them vanishes. This suggests that
$(A_f,\cL^{(u)})$ is not irreducible.
We may confirm this by giving its explicit decomposition.
Let
$$
\left(
\begin{matrix} \delta_1 \\ \delta_2\\ \delta_3\\ \delta_4
\end{matrix}
\right)=
\left(\begin{matrix}
1 & 0 & 0 & -1 \\
0 & 1 &-1 & 0\\
0 & 0 & 1 & 0 \\
0 & 0 & 0 & 1
\end{matrix}\right)
\left(
\begin{matrix} \gamma_1 \\ \gamma_2\\ \gamma_3\\ \gamma_4
\end{matrix}
\right)
$$
so that $(\delta_1, \delta_2, \delta_3, \delta_4)$ is a new symplectic basis for $E_u$.
The period matrix of $g_1, \overline{g_1}$ with respect to it is $\Omega'=(\Omega'_1,\mid \Omega'_2)$ with
$$
\begin{array}{c}
\Omega_1'=
\tiny{
\left(
\begin{matrix}
6.584352\dots + 8.916903\dots i &0  \\
0 & -15.444529\dots + 11.404432 \dots i
\end{matrix}
\right )  \,,} \\ \\
\Omega_2'=
\tiny{
 \left(
\begin{matrix}
-2.275825\dots + 6.429374\dots i &0 \\ \\
0 & -8.860177\dots + 2.487529\dots i
\end{matrix}
\right). }
\end{array}
$$
The lattices $\Lambda=\langle 6.584352 + 8.916903i, -2.275825 + 6.429374i  \rangle $,
$ \Lambda_{\sigma}=\langle -15.444529 + 11.404432i, -8.860177 + 2.487529i\rangle$ correspond to a pair of Galois
conjugate elliptic curves $C$ and ${}^{\sigma}C$ and the
 invariants of $C$ are\newline
$$
c_4(C)=\frac{9(5-12 \sqrt{-3})}{2^8}, \qquad
c_6(C)=\frac{27(43+42 \sqrt{-3})}{2^{12}},
$$
$$
j(C)=\frac{27(121171-36627 \sqrt{-3})}{686}.
$$
 Hence, we have seen that
 $$
 C\times \, ^{\sigma}C \stackrel{\Q(\sqrt{-3})}\simeq
  (A_f, \cL^{(u)})\stackrel{\Q}\simeq
\Res_{\Q(\sqrt{-3})/\Q}C.
$$

\subsection{Constructing a genus two curve from a quadratic $\Q$-curve}\label{constructing}

In the example above, we have modified the canonical polarization on the Jacobian of an hyperelliptic curve to present it as the Weil restriction of a $\Q$ curve. We now perform the reverse process, i.e., we depart from a quadratic $\Q$-curve and construct  a polarization on its Weil restriction which transforms it in the Jacobian of a
rational genus two curve. Notice that modular tools will not be used in this construction.

Let $C$ be the elliptic curve $Y^2 =X^3+aX+b$, where $a=-9(767+212 \sqrt{13})$ and
$b=-18(17225+4778\sqrt{13})$. We denote by $\sigma$ the nontrivial Galois conjugation of $K=\Q (\sqrt{13})$ over $\Q$. The points on $C$ with $x$-coordinate  equal to $3(-13+4\sqrt{13})$ generate a subgroup $G$ of 3-torsion points, and the elliptic curve $C/G$ is isomorphic over $K$ to  ${}^{\sigma}C$.
More precisely, there is a cyclic isogeny $\mu: C \rightarrow {}^{\sigma}C$ of degree $3$ defined over $K$ such that $\mu^{*}(\omega_{\sigma})=\lambda \omega$, where $\lambda=4+\sqrt{13}$ (see \cite{Go}) and $\omega_{\sigma}$
and $\omega$ are the invariant differential forms of ${}^{\sigma} C$ and $C$ respectively. In particular, ${}^{\sigma}\mu^{*}(\omega)={}^{\sigma}\lambda \omega_{\sigma}=3/\lambda\omega_{\sigma}$. We note that
the Weil restriction of $C$ is of real $\GL_2$-type, since  $\Q\otimes \End (\Res_{K/\Q} C)=\Q(\sqrt 3)$. Moreover, as there are no units of negative norm in this algebra, Corollary \ref{surpol} ensures the existence of a second polarization on $\Res_{K/\Q} C$. We will now build it.

Consider the period lattices
$$
\begin{array}{l}
\Lambda=\langle w_1=0.220377\dots ,w_2=0.428744 \dots i\rangle\,, \\
\Lambda_{\sigma}=\langle w_{\sigma,1}=-\lambda w_1 ,w_{\sigma,2}= -\lambda/3 w_2 \rangle \,,
\end{array}
$$
of the curves $C,\, ^{\sigma}C$ respectively.
Let $\gamma_1$, $\gamma_2$ (resp. $\gamma_{\sigma,1}$, $\gamma_{\sigma,2}$) be  a basis of $H_1(C,\Z)$
(resp. $H_1({}^{\sigma}C,\Z))$ such that $\int_{\gamma_i}\omega=w_i$
(resp. $\int_{\gamma_{\sigma,i}}\omega_{\sigma}=w_{\sigma,i}$). Then,   $(\gamma_1,\gamma_{\sigma,1}, \gamma_2,
\gamma_{\sigma,2})$ is a symplectic basis of
$H_1(C\times {}^{\sigma}C,\Z)$ for the canonical  polarization $\cL$ attached to $C\times {}^{\sigma}C$.

The action of the endomorphisms $T=\sqrt 3$ on $H_1(C\times {}^{\sigma}C,\Z)$ is obtained from its  action
 on $\Lambda \times \Lambda_{\sigma}$:
$$
\begin{array}{rcl}
T :\Lambda &\rightarrow &\Lambda_{\sigma}\\
w &\mapsto &
\lambda \cdot w\,,
\end{array}
\qquad
\begin{array}{rcl}
T: \Lambda_{\sigma} &\rightarrow&  \Lambda \\
 w &\mapsto &{}^{\sigma}
\lambda\cdot w=3/\lambda\cdot w\,.
\end{array}
$$
{}From this, we compute the matrix of the action of the
fundamental unit $u=2+\sqrt 3$ of $\End (\Res_{K/\Q} C)$ on
$H_1(C\times {}^{\sigma}C,\Z)$ (with respect to the basis
$\gamma_1,\gamma_{\sigma,1}, \gamma_2, \gamma_{\sigma,2}$):
$$M=\left (\begin{matrix}
 2  & -3 &  0  & 0\\
-1  &  2 &  0  & 0\\
 0  &  0 &  2  & -1\\
 0  &  0 & -3  &  2\end{matrix}
\right)\,.
$$
Hence, the Riemann form attached to the polarization $\cL^{(u)}$ is given by
$$E_u=\left (\begin{matrix}
 0 &  0 &  2 & -1 \\
 0 &  0 & -3 &  2 \\
-2 &  3 &  0 &  0 \\
 1 & -2 &  0 &  0
\end{matrix}
\right)\,.
$$
Taking into account that  $\omega_1=\omega+\omega_{\sigma}$, $\omega_2=(\omega-\omega_{\sigma})/\sqrt{13}$
form a rational basis of    $H^0(\Res_{K/\Q}C,\Omega^1_{\Q})$, we compute the period matrix of these differential
forms with respect to an $E_u$-symplectic basis of $H_1(C\times {}^{\sigma}C,\Z)$, and apply the procedure of \cite{GoGoGu}. We obtain that $(\Res_{K/\Q} C, \cL^{(u)})$ is the  Jacobian of one of the two curves
$$
Y^2=\pm (12909572X^6 + 17307966X^5 + 8746257X^4 + 2170636X^3
$$
$$ + 278850X^2 + 17238X+ 377
  )\,.
$$
Both curves have good reduction at $p=23$. Only for the curve corresponding to the sign $+$, we have that the characteristic polynomial of $\operatorname {Frob}_p$ acting on  the Tate module of  its Jacobian $\pmod p$ equals the square of the characteristic polynomial of $\operatorname {Frob}_p$ acting on the Tate module of $C/\F_p$. In conclusion, the right sign is $+$.

\subsection{Nonprincipally polarized abelian surfaces}

Coming back to the modular examples, we now illustrate what can be done to
describe briefly those abelian surfaces $A_f$ which are nonprincipally polarized.
While our ideas do not provide a systematic method to treat any abelian surface, since
we use the fact that they are of $\GL_2$-type, it covers many interesting cases. This will be apparent in the following subsection, where we will build a number of genus two curves with isogenous Jacobian.

We will work with the abelian surface $A_f$ which is the unique two-dimensional factor of $J_0(35)$. It corresponds to
the newform:
$$
f=q + aq^2 + (-a - 1)q^3 + (-a + 2)q^4 + q^5
- 4q^6 +\dots +(a+3)q^{13}+\dots
$$
where $a=(-1+\sqrt{17})/2$. We see that
 $\End_{\Q}(A_f)=\Z[ (1+\sqrt{17})/2]=:\cO$ is the integer ring of $K=\Q(a)$. Since
  $f$ does not have any extra-twist,
$A_f$ is simple in its $\overline{\Q}$-isogeny class and, in
particular, every principal polarization on $A_f$ determines a
curve whose Jacobian is isomorphic to  $A_f$. In addition, by
Corollary \ref{surpol}, if such a curve exists, it must be unique
up to $\qbar$-isomorphism.

Let $\sigma $ denote the Galois conjugation of $K/\Q$, and take
$\alpha=(17+\sqrt{17})/34$.
The cuspidal forms $h_1=\alpha f+\sigma(\alpha){}^{\sigma}f$, $h_2=(f-{}^{\sigma}f)/\sqrt{17}$
provide an integral basis of $H^0(A_f,\Omega_{\Q}^1)$.
A basis of  the $\Z$-module $H_1(A_f,\Z)$ is given by the  modular
symbols
$$
\begin{array}{l}
 \gamma_1=\{-1/10, 0\}  -\{-1/25, 0\}, \\
\gamma_2=    \{-1/21, 0\}  -\{-1/28, 0\},\\
\gamma_3=    \{-1/7, 0\} -\{-1/15, 0\} + \{2/5, 3/7\}, \\
\gamma_4=    \{-1/10, 0\}  -\{-1/15, 0\} + \{-1/25, 0\}  -\{-1/30, 0\}.
\end{array}
$$
In this basis, the matrix  of the
alternating Riemann form attached to the canonical polarization $\cL$ on $A_f$ is
$$
M_E=\left( E(\gamma_i,\gamma_j)\right)_{i,j}= \left(
\begin{matrix}
 0 &  1 & -1 & 0\\
-1 & 0 & 0 & 1\\
 1 & 0 & 0 & 1\\
 0 &-1 &-1 & 0\\
\end{matrix}\right).
$$
This polarization  over $\Q$ is not principal, since it is of type
$(1,2)$. Thus, we  cannot describe $(A_f,\cL)$ as a Jacobian or as a Weil restriction.
We will look for a different  polarization $\cL_0$ on $A_f$ allowing this explicit description for $A_f$.
Of course, we will require $\cL_0$ to be principal and defined over $\Q $.

The existence of this polarization is guaranteed by Proposition \ref{principal}.
By Theorem \ref{pi} and Corollary \ref{criterion}, we must check the polarizations $\cL^{(u)}$ for those endomorphisms $u\in \End_{\Q}(A_f)$  of norm 2. We take $u=T_{13}=(5+\sqrt{17})/2$. The rational representation of $u$
with respect to the basis $(\gamma_1,\gamma_2,\gamma_3,\gamma_4)$ is
$$
M_{u}= \left(
\begin{matrix}
 2 &1& 1& 1\\
 2 &2& 1& 2\\
 2 &0& 3&-2\\
 0 &1&-1& 3\\
 \end{matrix}
\right ).
$$
The symplectic product $E_u$ of the polarization $\cL ^{(u)}$ is
given by the matrix $M_{E_u}=M_E\cdot M_u$; it is of  type (2,2),
so that there is a principal polarization $\cL_0$ on $A_f$ over
$\Q$ such that $\cL_0^{\otimes 2}=\cL ^{(u)}$. A symplectic basis
for $H_1(A_f,\Z)$ with respect to this principal polarization is
$$
\left(
\begin{array}{c}
\delta_1 \\
\delta_2 \\
\delta_3 \\
\delta_4
\end{array}
\right)=\left(
\begin{array}{rrrr}
 1 & 0 & 0 & 0\\
 0 & 1 &1  &0\\
 0 & 2 &-1 &-1\\
 0 & 0 & 2 & 1
\end{array}\right)
\left(
\begin{array}{c}
\gamma_1 \\
\gamma_2 \\
\gamma_3 \\
\gamma_4
\end{array}
\right).
$$
We now compute the periods of the differential forms $h_1 dq/q, h_2 dq/q$ along
these paths to represent  $(A_f,\cL_0)$ as the
complex torus $\C /\Lambda $, where $\Lambda $ is the lattice
spanned by the columns of the matrix $\Omega=(\Omega_1\, |
\Omega_2)$, with
$$
\begin{array}{l}
\Omega_1=
\left(
\begin{matrix}
3.429722 \dots i & 1.265864 \dots i    \\ \\
- 0.224497 \dots i &  1.714858\dots i
\end{matrix}
\right),
\\ \\
\Omega_2=
\left(
\begin{matrix}
-4.737944 \dots&     3.044837\dots+1.898796\dots i    \\ \\
2.706904\dots &    -2.368972\dots+    2.572287\dots i
\end{matrix}
\right).
\end{array}
$$
We finally apply the method of \cite{GoGoGu} to see that this torus is the Jacobian of
the curve:
$$
\begin{array}{rl}
C_{35A}: Y^2=&
-(X+1)(8X+3)(10X^3 + 14X^2 + 6X + 1),
\end{array}
$$
whose absolute Igusa invariants are:
$$
\displaystyle \{i_1,i_2,i_3\}=\{-\frac{2^{23}29^5}{5^5 7^5},
\frac{2^{17} 29^3 37\cdot 83}{5^5 7^5},\frac{2^8 29^2 83^2 1913}{5^5 7^5}\}.
$$
As we mentioned before, this is the only rational curve whose Jacobian is isomorphic to $A_f$.

\subsection{Distinct genus two curves with isogenous Jacobians}

We now describe a method to find many nonisomorphic genus two curves with isogenous Jacobians.
As before, we will work with the only two-dimensional factor $A_f$ of $J_0(35)$.  We shall consider
abelian surfaces obtained as quotients $A_f/G$, where $G$ is a finite rational torsion subgroup
of $A_f$.

The existence of a nontrivial  rational 2-torsion point on the Jacobian
$J(C_{35,A})$ is evident from the equation of  the curve $C_{35,A}$. Rational
torsion points on $A_f$ can also be found by means of cuspidal
divisors on $X_0(35)$. We have the rational cusps $0,1/5,1/7,\infty$ of $X_0(35)$ at our
disposal.
Consider the cuspidal divisors $D_5=(0)-(1/5)$, $D_7=(0)-(1/7)$
and $D_{\infty}=(0)-(\infty)$ on $X_0(35)$. We denote by $G$
the group generated by their projections on $A_f$.

The integrals $\int_{\infty}^0 h_j dq/q$, provide the projection of $D_{\infty}$ on $A_f$: it is the point
$(-0.677587\dots, -0.084914\dots )$, which corresponds
to the {\it path} in $H_1(A_f,\Q)$ given by $1/8 \gamma_2-1/4 \gamma_3+1/8
\gamma_4= -9/8\delta_2+5/8\delta_3+3/4\delta_4$. Similarly, we
obtain
$$
\begin{array}{l}
D_5 \longleftrightarrow 1/8 \gamma_2+1/4 \gamma_3+1/8 \gamma_4= -1/8
\delta_2+1/8\delta_3+1/4\delta_4\,, \\
D_7\longleftrightarrow  1/2 \gamma_2=-1/2\delta_2+1/2\delta_3+1/2\delta_4\,.
\end{array}
$$
Thus $G\simeq \Z/8\Z \times \Z/2\Z\times \Z/2\Z$ is generated by the images of the divisors
$D_5,D_7, D_5-D_\infty$ (we shall denote these images by the same letters, since there is no risk of confusion). Note that all points in $G$ are rational torsion
points on $A_f$. For every subgroup $G'$ of $G$
the quotient abelian variety $A_f/G'$ could admit, in principle, a principal polarization
over $\Q$ since the degree of the polarization induced by $\cL$
on $A/G'$ is $2 \cdot \# G'\in \operatorname{Norm}(K_f/\Q)$. We will
only examine the cyclic subgroups of $G$.

The first step to check whether $A_f/G'$ is principally polarized
is the determination of
$\End_{\Q}(A_f/G')$. Let us consider
the Hecke operator $v=2u-1=T_{11}+T_{18}$
corresponding to the fundamental unit of negative norm $4+ \sqrt{17}$ in $K_f$. We
have:
$$
v (D_{\infty})=3 D_{\infty}\,, \quad v (D_5)=3 D_5\,,\quad v (D_7)=3 D_7\,,
$$
i. e., $v$ acts on $G$ as the multiplication by $3$ and hence leaves every subgroup of $G$
stable. This implies  that the order $\Z [\sqrt{17}]$  is contained in
$\End_{\Q}(A_/G')$ for all subgroups $G'$ of $G$. Nevertheless, the action of $u$
on $G$ only leaves the following cyclic subgroups stable:
\begin{center}
\begin{tabular}{|c|c|c|}
\hline
$G'=\langle P \rangle$ & order of $P$  &   $u(P)$\\
\hline
$\langle 4 D_{\infty}\rangle$ & $2$ &  $0$ \\
\hline
$\langle 2 D_{\infty}\rangle$ & $4$ &  $4 D_{\infty}$ \\
\hline
$ \langle D_5+D_7\rangle$ & $8$ & $-2 (D_5+ D_7)$  \\
\hline
\end{tabular}
\end{center}
For these three subgroups $G'$, we can ensure that $\End_{\Q}(A/G')=\cO$, and then
Proposition \ref{principal} tells us that $A/G'$ is
principally polarized; we will show how to build the principal polarization for   $A_f/\langle 4 D_{\infty}\rangle$.

In order to do so, we consider the lattice
$$
\Lambda'=\langle \gamma_1,\gamma_2,\gamma_3, \gamma_4,
1/2\gamma_2-\gamma_3+1/2\gamma_4\rangle
\subset H_1(A_f,\Q),
$$
with basis $\gamma'_1=1/2(\gamma_2+\gamma_4)$, $\gamma'_2=1/2(\gamma_2-\gamma_4)$,
$\gamma'_3=\gamma'_1$, $\gamma'_4=\gamma_3$.
The canonical polarization $E$ on $J_0(35)$ provides a natural symplectic form on
this lattice, which we shall denote by $E_{/\Lambda'}$. Its matrix with respect to the
basis $(\gamma_1', \gamma_2', \gamma_3', \gamma_4')$ is
$$M_{E_{/\Lambda'}}=\left(
\begin{matrix}
       0& -1/2& -1/2& -1/2\\
     1/2&    0& -1/2&  1/2\\
     1/2&  1/2&    0&   -1\\
     1/2& -1/2&    1&    0
\end{matrix}\right).
$$
The type of the corresponding primitive  polarization on
$A/\langle 4 D_{\infty}\rangle$ is $
(1,4)$. We now proceed as in the previous subsection to
derive a principal polarization from $E_{/\Lambda'}$: we use the  operator  $u^2$ (of norm 4)
to define the Riemann form $E'$ determined by  the matrix
$A/\langle 4 D_{\infty}\rangle$ given by the symplectic form
$$
M_{E'}:= M_{E_{/\Lambda'}}M_{u}^2=
\left(
\begin{matrix}
     0& -4 & -14& -14\\
     4& 0 & 6& 4\\
    14& -6 & 0& -8\\
     14&-4 & 8& 0
\end{matrix}\right),
$$
which is of type (2,2), and hence it is the square of a principal polarization $\cL'_0$.
Computing periods and jacobian Thetanullwerte, we find that
$(A_f/\langle 4 D_{\infty}\rangle, \cL'_0)$ is the Jacobian  of the curve:
$$
C_{35,B}: Y^2=-10X(4X+5)(5X+8)
    (25X^3 + 110X^2 + 156X + 70),
$$
with absolute Igusa invariants:
$$
\displaystyle \{i_1,i_2,i_3\}=\{
-\frac{2^{13} 43^5 359^5}{5^{17} 7^5},
-\frac{2^{10}43^3 359^3 8933}{5^{13} 7^5},$$ $$
-\frac{2^4 37\cdot43^2 359^2 571\cdot126949}{5^{13} 7^5}\}.
$$
Applying  the same procedure to the quotient
$A/\langle 2 D_{5}\rangle$, we arrive at the curve:
$$
C_{35,C}: Y^2=-2(11X + 16)(5X + 8)(3X + 5)(127X^3 + 598X^2 + 938X + 490),
$$
with absolute Igusa invariants:
$$
\displaystyle \{i_1,i_2,i_3\}=\{
\frac{2^{28}89^5}{5^{3} 7^7},
-\frac{2^{19}11\cdot23\cdot89^3 1489}{5^{3} 7^7},
-\frac{2^{10}43\cdot89^2 2683\cdot11239}{5^{3} 7^7}\}.
$$
Finally, the quotient $A/\langle  D_{5}+D_7\rangle$ is $\Q$-isomorphic to the
Jacobian of the
curve:
$$
\begin{array}{rl}
C_{35,D}: Y^2=\,& -142 (8X+13)(X^2 + 4744X + 3776)\\
& \qquad(2173X^3 + 10154X^2 + 15820X + 8218),
\end{array}
$$
with absolute Igusa invariants:
$$
\displaystyle \{i_1,i_2,i_3\}=\{
\frac{2^{13}109^5 214063^5}{5^{2} 7^5 71^{12}},
\frac{2^{12} 11\cdot17\cdot109^3 5171\cdot214063^3}{5^{2} 7^5 71^8},$$ $$
\frac{2^4 3^5 109^2 33871\cdot214063^2 271175273}{5^{2} 7^5 71^8}\}.
$$
In conclusion, we have found four nonisomorphic curves $C_{35,A}$, $C_{35,B}$, $C_{35,C}$, $C_{35,D}$, whose Jacobians are pairwise
isogenous.

\bibliographystyle{plain}
\bibliography{pola}

\end{document}